\documentclass[a4paper,fleqn]{cas-sc}

\usepackage{amsfonts}
\usepackage{amsopn}
\usepackage{graphicx}
\usepackage{epstopdf}
\usepackage{algorithm}
\usepackage{algorithmic}
\usepackage{cleveref}

\def\R{\mathbb{R}}
\def\Det{{\rm det}}
\def\Tr{{\rm Tr}}
\def\widebar#1{{\overline{#1}}}
\def\N{\mathbb{N}}

\newtheorem{theorem}{Theorem}[section]

\newtheorem{proposition}[theorem]{Proposition}
\newtheorem{remark}[theorem]{Remark}
\newtheorem{example}[theorem]{Example}
\newproof{proof}{Proof}

\crefname{proposition}{proposition}{propositions}
\Crefname{proposition}{Proposition}{Propositions}
\begin{document}
\let\WriteBookmarks\relax
\def\floatpagepagefraction{1}
\def\textpagefraction{.001}

\shorttitle{Efficient approximation of quadratic irrationals}
\shortauthors{van der Kamp et al.}

\title[mode=title]{On efficient approximation of quadratic irrationals, signed Chebyshev polynomials and Householder's method}

\author[1]{Peter H. {van der Kamp}}[orcid=0000-0002-2963-3528]
\cormark[1]
\ead{P.vanderKamp@LaTrobe.edu.au}
\affiliation[1]{organization={Dept. of Mathematical and Physical Sciences, La Trobe University},
            city={Victoria},
            postcode={3086},
            country={Australia}}

\author[1]{Anthony Overmars}[orcid=0000-0001-7021-4774]

\author[1]{Marcel Jackson}[orcid=0000-0002-8149-1141]

\author[2]{Andrew N.W. Hone}[orcid=0000-0001-9780-7369]
\affiliation[2]{organization={School of Engineering, Mathematics and Physics, University of Kent},
            city={Canterbury},
            country={UK}}

\cortext[1]{Corresponding author}

\begin{abstract}
We provide efficient algorithms to compute convergents of quadratic irrationals. We show that for square roots, in settings where Galois' refinement of Lagrange's theorem holds, certain decimations of the sequence of convergents are signed Chebyshev
sequences, which can be also be generated by a Householder method.
\end{abstract}


\maketitle


\section{Introduction}
Lagrange's theorem tells us that the simple continued fraction expansion of a quadratic irrational $\alpha$ is eventually periodic,
\begin{equation} \label{eq:e1}
\alpha=c_0+\frac{1}{c_1+\frac{1}{c_2+\cdots}}=[c_0,c_1,\ldots,c_r,\overline{c_{r+1},c_{r+2},\ldots,c_{r+l}}],
\end{equation}
i.e. $c_{r+i}=c_{r+i+l}$ for all $i>0$. The integer $l$ is the $l$ength of the period. This holds for the standard continued fraction for real numbers \cite[Theorem 4.3]{Olds}, but also for complex numbers \cite[Corollary 4.6]{DN}. In each case, the $m$th convergent of the continued fraction is
defined as
\begin{equation} \label{fpq}
\frac{p_m}{q_m}=[c_0,c_1,c_2,\ldots,c_m],
\end{equation}
and we have
\begin{equation} \label{eq:e2}
\begin{pmatrix}
c_0 & 1\\
1 & 0
\end{pmatrix}
\begin{pmatrix}
c_1 & 1\\
1 & 0
\end{pmatrix}
\cdots
\begin{pmatrix}
c_m & 1\\
1 & 0
\end{pmatrix}=
\begin{pmatrix}
p_m & p_{m-1}\\
q_m & q_{m-1}
\end{pmatrix},
\end{equation}
which implies, by taking determinants, that
\begin{equation} \label{eq:e3}
p_mq_{m-1}-p_{m-1}q_m=(-1)^{m+1},
\end{equation}
and hence each convergent \eqref{fpq} is in lowest terms. Convergents are the best rational approximations to $\alpha$ in the sense that there are no fractions closer to $\alpha$ with denominator $\leq q_m$.

In \Cref{sec:sec2} of this paper, we introduce Chebyshev polynomials, as well as dilated, sign-changed, and signed versions of those, and we prove some useful relations between them. 

In \Cref{sec:sec3}, we provide efficient algorithms to compute convergents of quadratic irrationals. Given the preperiod $r$, period $l$, and the initial period data, the $m$-th convergent, with
\[
m=r+\left(m-r \mod l\right)+l\left(2^{n_1}+2^{n_2}+\cdots+2^{n_q}\right),\quad n_1<n_2<\cdots<n_q,
\]
can be computed using $q+2$ matrix multiplications and $n_q$ matrix linear combinations, rather than generating the preceding $m-1$ convergents. This comes at the additional cost of calculating scalar coefficients, and it is compared with binary powering of the period matrix which amounts to $q+3+n_q$ matrix multiplications.

Galois showed that the continued fraction of a quadratic irrational $\alpha\in\R$, with conjugate $\overline{\alpha}$, is purely periodic if and only if $\alpha>1$ and $-1<\overline{\alpha}<0$ \cite[Theorems 4.1 and 4.2]{Olds}. As a consequence, the continued fraction of the square root of a rational number $N$ has the form \cite[Section 4.7]{Olds}, 
\begin{equation} \label{eq:sn}
\sqrt{N}=[c_0,\overline{c_1,c_2,\ldots,c_{2},c_{1},2c_0}]
\end{equation}
(note, if $N<1$ then $c_1$ plays the role of $c_0$). With $l$ the length of the period we have \cite[Equation (4.41)]{Olds}
\begin{equation} \label{eq:Pll}
p_{l-1}^2-Nq_{l-1}^2=(-1)^l.
\end{equation}
Galois' theorem has been extended to the setting of complex numbers in \cite{GGR}, where \eqref{eq:sn} was shown to hold, under additional conditions. In \Cref{sec:sec4}, we show that for continued fractions of the form \eqref{eq:sn}, the following decimations of the sequence of convergents,
\[
(p,q)_{kl-1}, \qquad k=0,1,\ldots,
\]
are signed Chebyshev sequences.

In \Cref{sec:sec5}, we show that these sequences of convergents can also be obtained by iterating a Householder method. 

\section{Dilated, sign-changed and signed Chebyshev polynomials}
\label{sec:sec2}
In this section, we introduce several sequences of polynomials related to the Chebyshev polynomials of the first and second kind. These can be defined in several ways: by linear recurrences, using trigonometric functions, as solutions of Pell equations, or via powers of two-by-two matrices with determinant $\pm1$. 
\begin{itemize}
\item The {\em Chebychev polynomials of the first and second kind}, $T_k$ and $U_k$, satisfy
\begin{align*}
T_{k+2}&=2xT_{k+1}-T_k,\qquad T_0=1,\ T_1=x\\
U_{k+2}&=2xU_{k+1}-U_k,\qquad U_0=1,\ U_1=2x.
\end{align*}
They can also be defined by
\begin{equation} \label{eq:dft}
T_k(\cos(\theta))=\cos(k\theta),\quad
U_k(\cos(\theta))=\frac{\sin((k+1)\theta)}{\sin(\theta)},
\end{equation}
or by
\[
e^{\iota k\alpha}=T_k(x)+U_{k-1}(x)y\iota, \quad e^{\iota\alpha}=x+\iota y,\qquad \iota^2=-1,
\]
or as solutions of the Pell equation
\[
T_k^2(x)-(x^2-1)U_{k-1}^2(x)=1.
\]
Of the many interesting properties we just mention the nesting properties
\begin{equation} \label{eq:nps}
T_{km}(x)=T_k(T_m(x)),\qquad U_{km-1}(x)=U_{k-1}(T_m(x))U_{m-1}(x),
\end{equation}
due to \eqref{eq:dft}. The former implies that the Chebyshev polynomials of the first kind, viewed as 1-dimensional maps, form a family of commuting maps, 
\[
T_a\circ T_b=T_b\circ T_a,
\]
and hence form an integrable system, see \cite{Ves}. 
\item The {\em dilated Chebychev polynomials of the first and second kind}\footnote{These polynomials were called modified Chebyshev polynomials of the first and second kind in \cite{RoDa} and Vieta-Lucas polynomials respectively  Vieta-Fibonacci polynomials, in \cite{Hor}.}, ${\cal T}_k$ and ${\cal U}_k$, see \cite{Andy}, satisfy
\begin{align*}
{\cal T}_{k+2}&=x{\cal T}_{k+1}-{\cal T}_k,\qquad {\cal T}_0=2,\ {\cal T}_1=x,\\
{\cal U}_{k+2}&=x{\cal U}_{k+1}-{\cal U}_k,\qquad {\cal U}_0=1,\ {\cal U}_1=x.
\end{align*}
They are monic polynomials, related to the Chebyshev polynomials by a scaling
\begin{equation} \label{eq:rdt}
{\cal T}_k(x)=2T_k\left(x/2\right),\quad {\cal U}_k(x)=U_k\left(x/2\right),
\end{equation}
and hence they possess the same nesting properties as the Chebyshev polynomials.
The dilated Chebyshev polynomial of the first kind arises naturally in \Cref{prop:pr1}.
\begin{proposition} \label[proposition]{prop:pr1}
Let $M$ be a $2\times 2$ matrix with determinant $\Det(M)=1$ and $\Tr(M)=x$. We have $\Tr\left(M^{k}\right)={\cal T}_k(x)$.
\end{proposition}
\begin{proof}
According to the Cayley--Hamilton theorem, $M$ satisfies its characteristic equation $M^2-\Tr(M)M+\Det(M){\bf 1}=0$, and hence
\[
\Tr\left(M^{k+2}\right)=\Tr\left(M^{k}(\Tr(M)M-\Det(M){\bf 1})\right)=x\Tr\left(M^{k+1}\right)-\Tr\left(M^{k}\right).
\]
Since $\Tr(M^0)=2$, the result follows.
\end{proof}
\item The {\em sign-changed dilated Chebyshev polynomials of the first and second kind}, $\widebar{T}_k(x)$ and $\widebar{U}_k$, can be obtained
from the dilated Chebyshev polynomials by changing each minus-sign in their coefficients to a plus-sign. They satisfy the sign-changed recursion
\begin{align*}
\widebar{{\cal T}}_{k+2}&=x\widebar{{\cal T}}_{k+1}+\widebar{{\cal T}}_k,\qquad \widebar{{\cal T}}_0=2,\ \widebar{{\cal T}}_1=x,\\
\widebar{{\cal U}}_{k+2}&=x\widebar{{\cal U}}_{k+1}+\widebar{{\cal U}}_k,\qquad \widebar{{\cal U}}_0=1,\ \widebar{{\cal U}}_1=x,
\end{align*}
and are related to the dilated Chebyshev polynomials by a complex scaling,
\begin{equation} \label{eq:rscdt}
\widebar{{\cal T}}_k(x)=\iota^k{\cal T}_k\left(x/\iota\right),\quad \widebar{{\cal U}}_k(x)=\iota^k{\cal U}_k\left(x/\iota\right).
\end{equation}
\begin{proposition} \label[proposition]{prop:pr2}
If $M$ is a $2\times 2$ matrix with determinant $\Det(M)=-1$ and $\Tr(M)=x$, then $\Tr\left(M^{k}\right)=\widebar{{\cal T}}_k(x)$.
\end{proposition}
The proof of \Cref{prop:pr2} is similar to the proof of \Cref{prop:pr1}.
\item The {\em sign-changed Chebyshev polynomials},
\begin{equation} \label{eq:rsct}
\widebar{T}_k=\iota^k{T}_k\left(x/\iota\right),\quad
\widebar{U}_k=\iota^k{U}_k\left(x/\iota\right),
\end{equation}
satisfy the sign-changed recursion
\begin{align*}
\widebar{T}_{k+2}&=2x\widebar{T}_{k+1}+\widebar{T}_k,\qquad \widebar{T}_0=1,\ \widebar{T}_1=x,\\
\widebar{U}_{k+2}&=2x\widebar{U}_{k+1}+\widebar{U}_k,\qquad \widebar{U}_0=1,\ \widebar{U}_1=2x.
\end{align*}
\item We define the {\em signed Chebyshev polynomials of the first and second kind},
\[
T_k^l(x)=\begin{cases}
T_k(x) & \text{when } l \equiv 0 \mod 2,\\
\widebar{T}_k(x) &\text{when } l \equiv 1 \mod 2,
\end{cases}\qquad
U_k^l(x)=\begin{cases}
U_k(x) &\text{when } l \equiv 0 \mod 2,\\
\widebar{U}_k(x) &\text{when } l \equiv 1 \mod 2,
\end{cases}
\]
which satisfy
\begin{align}
T^l_{k+2}&=2xT^l_{k+1}-(-1)^lT^l_k,\qquad T^l_0=1,\ T^l_1=x,\label{eq:Tl}\\
U^l_{k+2}&=2xU^l_{k+1}-(-1)^lU^l_k,\qquad U^l_0=1,\ U^l_1=2x.\label{eq:Ul}
\end{align}
as well as the {\em signed dilated Chebyshev polynomials of the first and second kind},
\[
{\cal T}_k^l(x)=\begin{cases}
{\cal T}_k(x) &\text{when } l \equiv 0 \mod 2,\\
\widebar{{\cal T}}_k(x) &\text{when } l \equiv 1 \mod 2,
\end{cases}\qquad
{\cal U}^l_k(x)=\begin{cases}
{\cal U}_k(x) &\text{when } l \equiv 0 \mod 2,\\
\widebar{{\cal U}}_k(x) &\text{when } l \equiv 1 \mod 2.
\end{cases}
\]
\item Finally, we introduce the Chebyshev polynomials of the third and fourth kind,
\begin{equation} \label{eq:dfw}
V_k(\cos(\theta))=\frac{\cos((2k+1)\theta/2)}{\cos(\theta/2)},\quad
W_k(\cos(\theta))=\frac{\sin((2k+1)\theta/2)}{\sin(\theta/2)},
\end{equation}     
which satisfy the recursion formulas
\begin{align*}
V_{k+2}&=2xV_{k+1}-V_k,\qquad V_0=1,\ V_1=2x-1,\\
W_{k+2}&=2xW_{k+1}-W_k,\qquad W_0=1,\ W_1=2x+1.
\end{align*}
\end{itemize}
The following relations will be useful.
\begin{proposition} \label[proposition]{prop:prp}
We have
\begin{enumerate}
\item $T_k(-T_2(x))=(-1)^kT_{2k}(x), \quad 2xU_k(-T_2(x))=(-1)^kU_{2k+1}(x)$,
\item $T_k(\widebar{T}_2(x))=\widebar{T}_{2k}(x), \quad 2xU_k(\widebar{T}_2(x))=\widebar{U}_{2k+1}(x)$,
\item $V_k(-T_2(x))=(-1)^kU_{2k}(x),\quad xW_k(-T_2(x))=(-1)^kT_{2k+1}(x)$,
\item $ V_k(\widebar{T}_2(x))=\widebar{U}_{2k}(x),\quad xW_k(\widebar{T}_2(x))=\widebar{T}_{2k+1}(x)$.
\end{enumerate}
\end{proposition}
\begin{proof}
Statement 1 readily follows from \cref{eq:nps} and the fact that the functions $T_k,U_k$ are even/odd when $k$ is even/odd. Combining formulas (14), (17) and (53) from \cite{Andy} one obtains
the relation
\[
V_k(x)=(-1)^kU_{2k}\left(\sqrt{\frac{1-x}2}\right)
\]
from which the first formula of statement 3 follows. The second formula of statement 3 is the case $l=1$ of the more general relation
\begin{equation} \label{eq:mgr}
T_lW_k(-T_{2l})=(-1)^kT_{(1+2k)l}
\end{equation}
which we prove using definitions \eqref{eq:dft} and \eqref{eq:dfw}, as follows. As $\cos(\phi)=-\cos(2l\theta)$ implies $\phi=\pi-2l\theta$, we have
\begin{align*}
T_lW_k(-T_{2l})(\cos(\theta))&=
\cos(l\theta)\frac{\sin((2k+1)(\pi-2l\theta)/2)}{\sin((\pi-2l\theta)/2)}\\
&=
\cos(l\theta)\frac{\cos((2k+1)l\theta-k\pi)}{\cos(l\theta)}\\
&=(-1)^k\cos((2k+1)l\theta).
\end{align*}
Statement 2 follows from statement 1, and statement 4 follows from statement 3, by using \eqref{eq:rsct}.
\end{proof}

\section{Efficient computation of convergents, when Lagrange holds} \label{sec:sec3}
Let $r,l$ be the lengths of the $r$un-up and the period respectively of a continued fraction \eqref{eq:e1}, and denote the matrix of convergents \eqref{eq:e2} by 
\[
\Psi_n:=\begin{pmatrix}
p_n & p_{n-1}\\
q_n & q_{n-1}
\end{pmatrix}.
\]
\begin{theorem} \label{thm:thm4}
We have, for all $i,k\in\mathbb{N}$ and $n\geq r$, 
\begin{equation}\label{eq:dobl}
\Psi_{n+2^{i}kl}=t_{2^{i-1}k}\Psi_{n+2^{i-1}kl}-(-1)^{2^{i-1}kl}\Psi_{n},
\end{equation}
where $t_j={\cal T}_j^l(t_1)$, in terms of
\begin{equation} \label{eq:et1}
t_1=(-1)^r\left(\left|
\begin{matrix}
p_{r+l-1} & p_r \\
q_{r+l-1} & q_r 
\end{matrix}\right|-\left|\begin{matrix}
p_{r+l} & p_{r-1} \\
q_{r+l} & q_{r-1} 
\end{matrix}\right|\right).
\end{equation}
\end{theorem}
\begin{proof}
From \eqref{eq:e1} and \eqref{eq:e2} we have, for all $n\geq r$,
\begin{align}
\Psi_{n+2l}&=\Psi_n\left(\Psi_n^{-1}\Psi_{n+l}\right)^2\notag\\
&=\Psi_n\left(\Tr(\Psi_n^{-1}\Psi_{n+l})\Psi_n^{-1}\Psi_{n+l}-\Det(\Psi_n^{-1}\Psi_{n+l})\right)\notag\\
&=\Tr(\Psi_n^{-1}\Psi_{n+l})\Psi_{n+l}-(-1)^l\Psi_{n},\label{eq:psir}
\end{align}
due to the Cayley--Hamilton theorem. We note that \cref{eq:psir} still holds if we replace $l$ by an integer multiple, so that for all $k\in\N$
\begin{equation}\label{eq:psirt}
\Psi_{n+2kl}=\Tr(\Psi_n^{-1}\Psi_{n+kl})\Psi_{n+kl}-(-1)^{kl}\Psi_{n}.
\end{equation}
The value
\begin{equation}\label{eq:dtk}
t_k:=\Tr\left(\Psi_n^{-1}\Psi_{n+kl}\right)=\Tr\left((\Psi_n^{-1}\Psi_{n+l})^k\right)
\end{equation}
does not depend on $n$, because, for all $t\geq r-n$, $\Psi_{n+t}^{-1}\Psi_{n+t+kl}$ is a conjugation of $\Psi_n^{-1}\Psi_{n+kl}$ which leaves the trace invariant. Since $\Det(\Psi_n^{-1}\Psi_{n+l})=(-1)^l$, it follows from \Cref{prop:pr1} and \Cref{prop:pr2} that
\begin{equation}\label{eq:tf}
t_k={\cal T}_k^l(t_1),
\end{equation}
in terms of the signed dilated Chebyshev polynomial of the first kind, i.e., we have
\begin{equation} \label{eq:tr}
t_{k+2}=xt_{k+1}-(-1)^lt_k,\quad t_0=2, t_1=x.
\end{equation}
From the definition \eqref{eq:dtk}, using \eqref{eq:e3}, we obtain the expression \eqref{eq:et1} for $t_1$.
Replacing $k\rightarrow 2^ik$ in \eqref{eq:psirt} yields \eqref{eq:dobl}.
\end{proof}
In the remaining of this section, we provide four algorithms to calculate the matrix of convergents $\Psi_m$, where
$m=m_0+lh$ with $m_0=r+\left(m-r \mod l\right)$, requiring $O(\log(h))$ operations.\footnote{Our count does not include the bit complexity caused by coefficient and matrix-entry growth}
The algorithms are based on the following base 2 representations of $h$:
\begin{itemize}
\item Binary. We write
\begin{equation} \label{eq:bin}
h=2^{n_1}+2^{n_2}+\cdots+2^{n_q},\quad n_1<n_2<\cdots<n_q
\end{equation}
and define
\begin{equation} \label{eq:hj}
h_j=m_0+l(2^{n_1}+\cdots+2^{n_j}),\quad j=1,\ldots,q,
\end{equation}
so that $h_q=m$.
\item Nested binary. We write
\begin{equation} \label{eq:ms}
h=2^{m_1}(1+2^{m_2}(1+\cdots+2^{m_{q-1}}(1+2^{m_q}))),
\end{equation}
with $m_1\geq 0$ and $m_i>0$ for $1<i\leq q$, and we define $k_0=1$,
\begin{equation} \label{eq:kj}
k_j=1+2^{m_{q+1-j}}k_{j-1},\quad j=1,\ldots,q-1.
\end{equation}
\end{itemize}
The numbers $n_i$ and $m_j$ are related by $n_i=m_1+m_2+\cdots+m_i$. 

The first two algorithms employ the linear relation \eqref{eq:dobl}. For each algorithm, the values $\Psi_{r+l}$, $\Psi_r$, $\Psi_{m_0}$, require precalculation directly by iteration.

\begin{algorithm}
\caption{Binary, additive.\\
Input: $r,l,m\geq r+l$, binary representation \eqref{eq:bin}, \eqref{eq:hj}.\\
Output: An expression for $\Psi_m$ in terms of $\Psi_{r+l}$, $\Psi_r$, $\Psi_{m_0}$ and the values \eqref{eq:valuesoft}.}
\label{alg:ALG1}
\begin{algorithmic}
\IF{$m_0\neq r$}
\STATE{Define $\Psi_{m_0+l}=\Psi_{r+l}\Psi_r^{-1}\Psi_{m_0}$}
\ENDIF
\STATE{Define $\Phi=\Psi_{m_0}^{-1}\Psi_{m_0+l}$}
\FOR{$j=1,2,\ldots,q$}
\FOR{$i=0,1,\ldots,n_{j-1}$}
\STATE{Define $\Psi_{h_{j-1}+l2^{i+1}}=t_{2^i}\Psi_{h_{j-1}+l2^i}-(-1)^{2^il}\Psi_{h_{j-1}}$}
\ENDFOR
\IF{$j\neq q$}
\STATE{Define $\Psi_{h_j+l}=\Psi_{h_j}\Phi$}
\ENDIF
\ENDFOR
\RETURN $\Psi_m$
\end{algorithmic}
\end{algorithm}

\begin{proposition}
The cost of \Cref{alg:ALG1} is $q+2$ (or $q$ if $m_0=r$) matrix multiplications, $\sum_{i=1}^q n_i$ linear combinations of matrices, and one has to determine the values
\begin{equation} \label{eq:valuesoft}
t_{2^i},\quad i=0,1,\ldots,n_q-1.
\end{equation}
\end{proposition}
The values \eqref{eq:valuesoft} can be determined using
\begin{align}
t_{2^{i+1}}&=\Tr\left(((\Psi_n^{-1}\Psi_{n+l})^{2^i})^2\right)\notag\\
&=t_{2^i}^2-(-1)^{2^il}2.\label{t2}
\end{align}

\begin{example} \label{exa:ex1}
Consider
$\alpha=\frac{4}{3}+\frac{\sqrt{3}}{6}=[1, 1, 1, 1, \overline{1, 1, 4, 1, 1, 2, 20, 2}]$.
So we have $r=3$, $l=8$. We show how to calculate $\Psi_m$ for $m=89=9+8(2+2^3)$. We use \eqref{eq:e2} to calculate $\Psi_i$ for $i=0,1,\ldots,l+r$. In particular we have, for $i=r,m_0,l+r$
\[
\Psi_3=\left[\begin{array}{cc}
5 & 3 
\\
 3 & 2 
\end{array}\right], \quad
\Psi_9=\left[\begin{array}{cc}
339 & 133 
\\
 209 & 82 
\end{array}\right]
, \quad
\Psi_{11}=\left[\begin{array}{cc}
14165 & 6913 
\\
 8733 & 4262 
\end{array}\right],
\]
and hence $t_1=\Tr(\Psi_3^{-1}\Psi_{11})=2702$. We then use \eqref{t2} to determine
\begin{align*}
t_{2} &= t_1^2-2 = 7300802,\\
t_{4} &= t_2^2-2 = 53301709843202.
\end{align*}
We subsequently calculate
\begin{align*}
\Psi_{17}&=\Psi_{11}\Psi_3^{-1}\Psi_9,\\
\Phi&=\Psi_9^{-1}\Psi_{17},\\
\Psi_{25}&=t_1\Psi_{17}-\Psi_9,\\
\Psi_{33}&=\Psi_{25}\Phi,\\
\Psi_{41}&=t_1\Psi_{33}-\Psi_{25},\\
\Psi_{57}&=t_2\Psi_{41}-\Psi_{25},\\
\Psi_{89}&=t_4\Psi_{57}-\Psi_{25},
\end{align*}
to find
\[
\Psi_{89}={\small \left[\begin{array}{cc}
7031582616783360742995441537263465239 & 2758523931487789014011972217814706733 
\\
 4335108450922621626554341085216343809 & 1700684050688932407684112398936807682 
\end{array}\right]}.
\]
\end{example}

Next we give a Horner-type algorithm. 
\begin{algorithm}
\caption{\label{alg:ALG2} Nested binary, additive.\\
Input: $r,l,m\geq r+l$ and nested binary representation \eqref{eq:ms},\eqref{eq:kj}.\\
Output: An expression for $\Psi_m$ in terms of $\Psi_{r+l}$, $\Psi_r$, $\Psi_{m_0}$ and the values \eqref{eq:vot}.}

\begin{algorithmic}
\IF{$m_0\neq r$}
\STATE{Define $\Psi_{m_0+l}=\Psi_{r+l}\Psi_r^{-1}\Psi_{m_0}$}
\ENDIF
\STATE{Define $\Phi=\Psi_{m_0}^{-1}\Psi_{m_0+l}$}
\FOR{$j=0,\ldots,q-1$}
\FOR{$i=0,\ldots,m_{q-j}-1$}
\STATE{Define $\Psi_{m_0+2^{i+1}k_jl}=t_{2^ik_j}\Psi_{m_0+2^ik_jl}-(-1)^{2^ik_jl}\Psi_{m_0}$}
\ENDFOR
\IF{$j\neq q-1$}
\STATE{Define $\Psi_{m_0+lk_{j+1}}=\Psi_{m_0+2^{m_{q-j}}lk_j}\Phi$}
\ENDIF
\ENDFOR
\RETURN $\Psi_m$
\end{algorithmic}
\end{algorithm}

\begin{proposition}
The cost of \Cref{alg:ALG2} is $q+2$ (or $q$ if $m_0=r$) matrix multiplications and $\sum_{i=1}^q m_i=n_q$ linear combinations of matrices, and one has to determine the values
\begin{equation} \label{eq:vot}
t_{k_{j-1}2^i},\quad i=0,1,\ldots,m_{q+1-j}-1,\quad j=1,\ldots,q.
\end{equation}
\end{proposition}
The cost is less than the $\sum_{i=1}^q n_i=\sum_{i=1}^q (q+1-i)m_i$ that were required in \Cref{alg:ALG1}. The values \eqref{eq:vot} can be determined using \Cref{alg:ALG3}, which is based on the formula
\begin{equation} \label{nice}
\begin{pmatrix}
t_{2k}\\
t_{2k+1}
\end{pmatrix}=t_k
\begin{pmatrix}
t_{k}\\
t_{k+1}
\end{pmatrix}
-(-1)^{kl}\begin{pmatrix}
t_{0}\\
t_{1}
\end{pmatrix}.
\end{equation}

\begin{algorithm}
\caption{\label{alg:ALG3} Determining coefficients.\\
Input: $q,m_1,\ldots,m_q,k_0,\ldots,k_{q-1}$.\\
Output: The values \eqref{eq:vot}.}
\begin{algorithmic}
\IF{$m_1=0$}
\STATE{Define $p=q-1$}
\ELSE
\STATE{Define $p=q$}
\ENDIF
\FOR{$j=1,\ldots,p$}
\IF{$j=1$}
\STATE{Define $t_2=t_1^2-(-1)^l2$}
\ELSE
\STATE{Define $z=(k_{j-1}-1)/2$}
\IF{$z=1$ {\bf or} $j=p$}
\STATE{Define $t_{2z+1}=t_zt_{z+1}-(-1)^{lz}t_1$}
\ELSE 
\STATE{Define $\begin{aligned}
\begin{pmatrix}
t_{2z}\\
t_{2z+1}
\end{pmatrix}&=t_z
\begin{pmatrix}
t_{z}\\
t_{z+1}
\end{pmatrix}
-(-1)^{lz}\begin{pmatrix}
2\\
t_{1}
\end{pmatrix}\\
t_{2z+2}&=t_1t_{2z+1}-(-1)^lt_z
\end{aligned}$}
\ENDIF  
\ENDIF
\FOR{$i=1,\ldots,m_{q+1-j}-1$}
\STATE{Define $z=k_{j-1}2^{i-1}$}
\IF{$z=1$}
\STATE{Define $t_{2z+1}=t_zt_{z+1}-(-1)^{lz}t_1$}
\ELSE
\IF{$j=p$}
\STATE{Define $t_{2z}=t_z^2-(-1)^{lz}2$}
\ELSE 
\STATE{Define $
\begin{pmatrix}
t_{2z}\\
t_{2z+1}
\end{pmatrix}=t_z
\begin{pmatrix}
t_{z}\\
t_{z+1}
\end{pmatrix}
-(-1)^{lz}\begin{pmatrix}
2\\
t_{1}
\end{pmatrix}
$}
\ENDIF
\ENDIF
\ENDFOR
\ENDFOR
\RETURN $t_{k_{j-1}2^i},\quad i=0,1,\ldots,m_{q+1-j}-1,\quad j=1,\ldots,q$.
\end{algorithmic}
\end{algorithm}

\begin{example}
For the same $\alpha,m$ as in \Cref{exa:ex1}, we show how to find $\Psi_m$ using \Cref{alg:ALG2} instead of \Cref{alg:ALG1}. We write $m=89=9+2(1+2^2)8$ so that $q=2$, $m_0=9$, $m_1=1,m_2=2$ and $k_0=1,k_1=5$. According to \cref{eq:vot} we need to determine $t_1,t_2,t_5$. Using \Cref{alg:ALG3} this is done as follows:
\begin{align*}
t_{2} &= t_{1}^{2}-t_{0} = 7300802 \\
t_{3} &= t_{1} t_{2}-t_{1} = 19726764302 \\
t_{5} &= t_{2} t_{3}-t_{1} = 144021200269567502.
\end{align*}
After the initial calculation of $\Psi_3,\Psi_9,\Psi_{11}$, \Cref{alg:ALG2} determines
\begin{align*}
\Psi_{17}&=\Psi_{11}\Psi_3^{-1}\Psi_9,\\ \Phi&=\Psi_9^{-1}\Psi_{17},\\
\Psi_{25}&=t_1\Psi_{17}-\Psi_9,\\
\Psi_{41}&=t_2\Psi_{25}-\Psi_9,\\
\Psi_{49}&=\Psi_{41}\Phi,\\
\Psi_{89}&=t_5\Psi_{49}-\Psi_9.
\end{align*}
\end{example}

The Hurwitz continued fraction for a complex number $\alpha$ is obtained by setting $\alpha_0=\alpha$ and
\begin{equation} \label{Alg}
c_i=\left[ \alpha_i \right], \quad \alpha_{i+1}=\left(\alpha_i-c_i \right)^{-1}, \quad i=0,1,\ldots,
\end{equation}
where $\left[ \alpha \right]$ denotes rounding $\alpha$ to the nearest Gaussian integer (up in case of a tie). Galois' theorem was proven to hold under certain constraints in \cite{GGR}. As a consequence, the Hurwitz continued fraction of a square root often (but not always) takes the form \eqref{eq:sn}.

\begin{example}
The Hurwitz continued fraction of $\sqrt{3+10\iota}$ is
\[
[3+2 \iota, \overline{-2, -1-\iota, -1+2 \iota, -2-\iota, 2, 5+4 \iota}]
\].
We observe that $c_6\neq 2c_0$, $c_5\neq c_1$ and $c_4\neq c_2$. One can still apply any of the algorithms \ref{alg:ALG1}, \ref{alg:ALG2}, \ref{alg:ALG3} or \ref{alg:ALG4}.
\end{example}

\begin{example}
The Hurwitz continued fraction of $\sqrt{9+10\iota}$ is
\begin{align*}
\sqrt{9+10\iota}=[3 + \iota,\overline{1 - \iota, 3\iota, -2 + \iota, -2 - \iota, 3 - 2\iota, -2 + 3\iota,
3 - 2\iota, -2 - \iota, -2 + \iota, 3\iota, 1 - \iota, 6 + 2\iota}&].
\end{align*}
It has the form \eqref{eq:sn} with period length 12. We use \Cref{alg:ALG2} to determine $\Psi_{71}$. Writing 71 in the form \eqref{eq:ms} gives $q=2$, $m_0=11,m_1=0,m_2=2$ and $k_0=1,k_1=5$. We have
\begin{equation} \label{eq:psis}
\Psi_0=\begin{pmatrix}
3+\iota  & 1 
\\
 1 & 0 
\end{pmatrix},\quad
\Psi_{11}=\begin{pmatrix}
-101025+51393 \iota  & -60722-31709 \iota  
\\
 -19460+24005 \iota  & -18640-1162 \iota  
\end{pmatrix}.
\end{equation}
We find $t_1=\Tr(\Psi_0^{-1}\Psi_{12})=-202050 + 102786\iota$ and
$t_2=t_1^2-2=30259240702 - 41535822600\iota$.
Following \Cref{alg:ALG2} we calculate
\begin{align*}
\psi_{23} &= \psi_{12}\psi_{0}^{-1} \psi_{11}, \\
\Phi &= \psi_{11}^{-1} \psi_{23}, \\
\psi_{35} &= t_{1} \psi_{23}-\psi_{11}, \\
\psi_{59} &= t_{2} \psi_{35}-\psi_{11}, \\
\psi_{71} &= \psi_{59}\Phi,
\end{align*}
from which we extract
\begin{equation} \label{eq:pq71}
{\small
\begin{split}
p_{71}&=-64452969879034582258134562726849-
21217336886334890599158733121700 \iota,\\
q_{71}&=-18405487633517442616165619582790+1864795250277698166333066426570 \iota.
\end{split}}
\end{equation}
\end{example}

Both \Cref{alg:ALG1} and \Cref{alg:ALG2} employ the linear relation \eqref{eq:dobl} from \Cref{thm:thm4} as much as possible, minimising the use of matrix multiplication (by $\Phi$). The next two algorithms are solely based on matrix multiplication.
\begin{proposition}
The cost of Algorithms \ref{alg:ALG4} and \ref{alg:ALG5}, is $n_q+q+3$ matrix multiplications (or $n_q+q+1$ if $m_0=r$). (And one does not need to determine any coefficients.)
\end{proposition}

\begin{algorithm}
\caption{Binary, multiplicative.\\
Input: $r,l,m\geq r+l$, binary representation \eqref{eq:bin}, \eqref{eq:hj}.\\
Output: An expression for $\Psi_m$ in terms of $\Psi_{r+l}$, $\Psi_r$, $\Psi_{m_0}$.}
\label{alg:ALG4}
\begin{algorithmic}
\IF{$m_0\neq r$}
\STATE{Define $\Psi_{m_0+l}=\Psi_{r+l}\Psi_r^{-1}\Psi_{m_0}$}
\ENDIF
\STATE{Define $\Phi=\Psi_{m_0}^{-1}\Psi_{m_0+l}$, and $n_0=0$.}
\FOR{$j=1,2,\ldots,q$}
\FOR{$i=n_{j-1}+1,n_{j-1}+2,\ldots,n_{j}$}
\STATE{Redefine $\Phi\leftarrow \Phi^2$}
\ENDFOR
\STATE{Define $\Psi_{h_{j}}=\Psi_{h_{j-1}}\Phi$.}
\ENDFOR
\RETURN $\Psi_m$.
\end{algorithmic}
\end{algorithm}

\begin{algorithm}
\caption{Nested binary, multiplicative.\\
Input: $r,l,m\geq r+l$, nested binary representation \eqref{eq:ms}, \eqref{eq:kj}.\\
Output: An expression for $\Psi_m$ in terms of $\Psi_{r+l}$, $\Psi_r$, $\Psi_{m_0}$.}
\label{alg:ALG5}
\begin{algorithmic}
\IF{$m_0\neq r$}
\STATE{Define $\Psi_{m_0+l}=\Psi_{r+l}\Psi_r^{-1}\Psi_{m_0}$}
\ENDIF
\STATE{Define $\Theta=\Phi=\Psi_{m_0}^{-1}\Psi_{m_0+l}$}
\FOR{$j=q,q-1,\ldots,1$}
\FOR{$i=1,2,\ldots,m_{j}$}
\STATE{Redefine $\Theta\leftarrow \Theta^2$}
\ENDFOR
\IF{$j\neq1$}
\STATE{Redefine $\Theta\leftarrow \Theta\Phi$}
\ENDIF
\ENDFOR
\STATE{Define $\Psi_m=\Psi_{m_0}\Theta$}
\RETURN $\Psi_m$.
\end{algorithmic}
\end{algorithm}

\Cref{alg:ALG1} is the most direct binary method, but in the worst case its $\sum_{i=1}^q n_i$ matrix linear combinations give quadratic growth in $\log h$. \Cref{alg:ALG2} avoids this by using the nested form; it is the preferred method. \Cref{alg:ALG2} will often be faster than the multiplicative variants, because it replaces the $n_q$ squarings of a $2\times2$ matrix by $n_q$ linear combinations of such matrices.

\section{When Galois holds, a decimation of the sequence of convergents is a Chebyshev sequence}
\label{sec:sec4}
Throughout this section, and the next, we consider numbers $\alpha=\sqrt{N}$ for which the continued fraction has the form~\eqref{eq:sn}. For such numbers, certain decimations (a decimation is an arithmetic subsequence) of the sequence of convergents are signed Chebyshev sequences.
\begin{theorem} \label{thm:tcs}
Let $p_i,q_i$, $i=0,1,\ldots$, be the numerators and denominators of the convergents of a continued fraction of the form \eqref{eq:sn}, so that equation \eqref{eq:Pll} holds. We then have
\begin{equation} \label{eq:pqC}
\begin{split}
p_{kl-1}&=T^l_k(p_{l-1}),\\
q_{kl-1}&=q_{l-1}U^l_{k-1}(p_{l-1}).
\end{split}
\end{equation}
\end{theorem}
\begin{proof}
The palindromic nature of \eqref{eq:sn} tells us that the matrix $\Psi_0^{-1}\Psi_{l-1}$ is symmetric, which yields $q_{l-2}=p_{l-1}-c_0q_{l-1}$. Together with the fact that
$c_l=2c_0$, this implies that
\begin{align*}
t_1&=\Tr\left(\Psi_0^{-1}\Psi_{l}\right)\\
&=\Tr\left(
\begin{pmatrix} 0 & 1 \\ 1 & -c_0 \end{pmatrix}
\begin{pmatrix} p_{l-1} & p_{l-2} \\ q_{l-1} & q_{l-2} \end{pmatrix}
\begin{pmatrix} 2c_0 & 1 \\ 1 & 0 \end{pmatrix}\right)\\
&=c_0q_{l-1}+p_{l-1}+q_{l-2}\\
&=2p_{l-1}.
\end{align*}
The ${}_{1,1}$-component of formula \eqref{eq:psir} tells us that, with $p_{l-1}=x$, for all $n$,
\[
p_{n+2l}=2xp_{n+l}-(-1)^lp_{n}.
\]
Since $p_{-1}=1$, according to \eqref{eq:Tl} we have $p_{kl-1}=T^l_k(x)$.
Similarly, the ${}_{2,1}$-component of formula \eqref{eq:psir} yields
\[
q_{n+2l}=2xq_{n+l}-(-1)^lq_{n},
\]
and, as $q_{-1}=0$, this implies $q_{2l-1}=2xq_{l-1}$. Thus, we see that up to a constant factor $q_{l-1}$, the decimation $q_{kl-1}$ is described by a signed Chebyshev polynomial of the second kind, $U^l_{k-1}(x)$, cf. \eqref{eq:Ul}.
\end{proof}

\begin{example}
We verify that the 71$^{st}$ convergent of $\sqrt{9+10\iota}$ is
given by equation \eqref{eq:pqC}. From \eqref{eq:psis}, we obtain
\[
p_{11}=-101025+51393 \iota,\quad  q_{11}= -19460+24005 \iota.
\]
Substituting $x=p_{11}$ into $T_6$ and $q_{11}U_5$, where
\[
T_6=(2x^2 - 1)(16x^4 - 16x^2 + 1), \quad U_5=2x(2x + 1)(2x - 1)(4x^2 - 3),
\]
yields \eqref{eq:pq71}.
\end{example}

\section{Householder's method} \label{sec:sec5}
\begin{theorem} \label{thm:HHT}
Let $f(x)=x^2-N$ and denote by $(1/f)^{(d)}$ the $d$-th derivative of $1/f$. The order $d$ numerical method from Householder,
\begin{equation} \label{eq:HHE}
H(x)=x+d\frac{(1/f)^{(d-1)}(x)}{(1/f)^{(d)}(x)},
\end{equation}
applied to the $(l-1)$st convergent of $\sqrt{N}$, $x=p_{l-1}/q_{l-1}$, yields the $(kl-1)$st convergent, where $k=d+1$.
\end{theorem}
\begin{proof}
According to \cite[Theorem 2]{YD}, we have, with $X=\frac{N+x^2}{N-x^2}$, 
\begin{equation} \label{eq:Hx}
H(x)=\begin{cases}
x\frac{W_{d/2}(X)}{V_{d/2}(X)} & \text{when } d \text{ is even},\\[1mm]
x\frac{T_{(d+1)/2}(X)}{(X-1)U_{(d-1)/2}(X)}
& \text{when } d \text{ is odd}.
\end{cases}
\end{equation}
For convenience we write $x=p/q$. We have, using \eqref{eq:Pll},
\[
X=\frac{Nq^2+p^2}{Nq^2-p^2}=1-(-1)^l2p^2=
\begin{cases}
-T_2(p) & \text{when } l \text{ is even},\\ 
\widebar{T}_2(p) & \text{when } l \text{ is odd}.
\end{cases}
\]
Substituting the appropriate expression for $X$ in \eqref{eq:Hx}, and using statement $i$ of \Cref{prop:prp} with
\[
i=\begin{cases}
1 & \text{when } d \text{ is odd and } l \text{ is even},\\
2 & \text{when } d \text{ is odd and } l \text{ is odd},\\
3 & \text{when } d \text{ is even and } l \text{ is even},\\ 
4 & \text{when } d \text{ is even and } l \text{ is odd},\\ 
\end{cases}
\]
yields
\[
H(x)=\frac{T^l_{d+1}(p)}{qU^l_d(p)},
\]
which according to \Cref{thm:tcs} equals $p_{kl-1}/q_{kl-1}$ with $k=d+1$.
\end{proof}

\begin{remark}
The special case $k=2$ of \eqref{eq:pqC}, which corresponds to $d=1$ of \eqref{eq:HHE}, is called Newton's method,
\begin{equation} \label{Newton}
\frac{p_{2l-1}}{q_{2l-1}}=\frac{T^l_2(p_{l-1})}{q_{l-1}U^l_{1}(p_{l-1})}
=\frac{2p_{l-1}^2-(-1)^l}{2q_{l-1}p_{l-1}}.
\end{equation}
It coincides with the Babylonian (or Heron's) method, cf. \cite{YD} and \cite{O1}, where the method was used to approximate $\sqrt{2}=[1,\overline{2}]$, i.e., one iterates
\[
\frac{p}{q}\mapsto\frac{1}{2}\frac{p}{q}+\frac{q}{p}.
\]
Note that the right hand side equals
\[
\frac{2p^{2}-1}{2 q p}=\frac{T_2(p)}{qU_{1}(p)},
\]
if $p/q$ satisfies the Pell equation $p^2-2q^2=1$. The case $k=3$ is Halley's method,
\begin{equation} \label{Halley}
\frac{p_{3l-1}}{q_{3l-1}}=\frac{T^l_3(p_{l-1})}{q_{l-1}U^l_{2}(p_{l-1})}
=\frac{p_{l -1}}{q_{l -1}} \frac{4 p_{l -1}^{2}-3 \left(-1\right)^{l}}{4p_{l-1}^{2}-\left(-1\right)^{l}}.
\end{equation}
\end{remark}

\begin{remark}
For large $k$, equation \eqref{eq:pqC} can be evaluated in $O(\log k)$ matrix multiplications, see \cite[Section 3.1]{CLH}, by writing the signed Chebyshev recursion in the form $\Xi^l_n=M_l\Xi^l_{n-1}$, where
\begin{equation} \label{eq:mbm}
\Xi^l_n=
\begin{pmatrix}
T^l_{n+1} & U^l_{n+1} \\
T^l_{n} & U^l_{n}
\end{pmatrix},\quad
M_l=
\begin{pmatrix}
2x & -(-1)^l \\
1 & 0
\end{pmatrix},\quad
\Xi^l_0=\begin{pmatrix}
x & 2x \\
1 & 1
\end{pmatrix}.
\end{equation}
Alternatively, one may apply the halve-and-square method, see \cite[Section 3.2]{CLH}, based on
\begin{equation} \label{has}
\Xi^l_{n}=\begin{cases}
2T^l_m\Xi^l_m-(-1)^{ml}\Xi^l_0 & n = 2m, \\
2T^l_m\Xi^l_{m+1}-(-1)^{ml}\Xi^l_1 & n = 2m +1, \\
\end{cases}
\end{equation}
which follows from the Cayley--Hamilton theorem and the facts that Tr$(M_l^m)=2T^l_m$ and Det$(M_l^m)=(-1)^{ml}$. 
\end{remark}

\appendix
\section*{Acknowledgments}
The first and last author are grateful for the hospitality of Da-jun Zhang and Cheng Zhang from Shanghai University and J.P. Wang from Ningbo University, where part of this research was performed. P.H. van der Kamp was supported by NSFC grant No. 12271334.


\begin{thebibliography}{10}
\bibitem{Andy}
A.N.W. Hone, L.E. Jeffery, R.G. Selcoe,
On a family of sequences related to Chebyshev polynomials,
J. Integer Seq. 21 (2018) 18.7.2.

\bibitem{CLH}
Z. Cao, L. Liu, L, Hong,
Evaluation methods for Chebyshev polynomials,
Cryptology {ePrint} Archive 2020/1365
(2020), https://eprint.iacr.org/2020/1365.

\bibitem{DN}
S.G. Dani and A. Nogueira, Continued fractions for complex numbers
and values of values of binary quadratic forms, Trans. Amer. Math. Soc. 366(7) (2014) 3553--3583.

\bibitem{YD}
Y. Dijoux, Chebyshev polynomials involved in the Householder’s method for square roots, Carpathian Math. Publ. 17 (2025) 616--630.

\bibitem{GGR}
G. Gonz\'alez Robert, Purely periodic and transcendental complex continued fractions,
Acta Arith. 194 (2020) 241--265.

\bibitem{Hor}
A.F. Horadam, Vieta polynomials, Fib. Quart. 40 (2002) 223--232.

\bibitem{OEIS}
N. J. A. Sloane, ed., The On-Line Encyclopedia of Integer Sequences, published electronically at https://oeis.org/A010337 (accessed 15 Nov. 2025).

\bibitem{Olds}
C.D. Olds, Continued Fractions, New mathematical library 9(4), Random House, 1963.

\bibitem{O1}
A. Overmars, S. Venkatraman, S. Parvin, Revisiting square roots with a fast estimator, Lond. J. Res. Comput. Sci. Technol. 18(1) (2018) 45--52.

\bibitem{RoDa}
R. Witu\l a and D. S\l ota, On modified Chebyshev polynomials, J. Math. Anal. Appl. 324 (2006) 321--343.

\bibitem{Ves}
A.P. Veselov, Integrable maps, Russ. Math. Surveys 46 (1991) 1--51.
\end{thebibliography}
\end{document}